\RequirePackage{ifpdf}
\ifpdf 
\documentclass[pdftex]{sigma}
\else
\documentclass{sigma}
\fi

\begin{document}
\allowdisplaybreaks

\renewcommand{\PaperNumber}{067}

\FirstPageHeading

\ShortArticleName{The Relation Between the Associate Almost
Complex Structure}

\ArticleName{The Relation Between the Associate Almost Complex\\
Structure to $\boldsymbol{HM'}$ and
$\boldsymbol{(HM',S,T)}$-Cartan\\ Connections}

\Author{Ebrahim ESRAFILIAN and Hamid Reza SALIMI MOGHADDAM}
\AuthorNameForHeading{E. Esraf\/ilian and H.R. Salimi Moghaddam}

\Address{Department of Pure Mathematics, Faculty of Mathematics,\\
Iran University of Science and Technology, Narmak-16, Tehran,
Iran}

\Email{\href{mailto:salimi_m@iust.ac.ir}{salimi\_m@iust.ac.ir}}

\ArticleDates{Received April 08, 2006, in f\/inal form August 30,
2006; Published online September 06, 2006}

\Abstract{In the present paper, the $(HM',S,T)$-Cartan connections
on pseudo-Finsler manifolds, introduced by A.~Bejancu and
H.R.~Farran, are obtained by the natural almost complex structure
arising from the nonlinear connection $HM'$. We prove that the
natural almost complex linear connection associated to a
$(HM',S,T)$-Cartan connection is a metric linear connection with
respect to the Sasaki metric $G$. Finally we give some conditions
for $(M', J, G)$ to be a K\"ahler manifold.}

\Keywords{almost complex structure; K\"ahler and pseudo-Finsler
manifolds; $(HM',S,T)$-Cartan connection}

\Classification{53C07; 53C15; 53C60; 58B20}

\section{Introduction}

Almost complex structures are important structures in dif\/ferential
geometry~\cite{[Ic1],[Ic2],[Ma]}. These structures have found many
applications in physics. H.E.~Brandt has shown that the spacetime
tangent bundle, in the case of Finsler spacetime manifold, is
almost complex~\cite{[Brandt2],[Brandt3],[Brandt4]}. Also he
demonstrated that in this case the spacetime tangent bundle is
complex provided that the gauge curvature f\/ield
vanishes~\cite{[Brandt1]}. In \cite{[BeFa], [BeFa1]}, for a
pseudo-Finsler manifold $F^m=(M,M',F^{\ast})$ with a nonlinear
connection $HM'$ and any two skew-symmetric Finsler tensor f\/ields
of type $(1,2)$ on $F^m$, A.~Bejancu and H.R.~Farran introduced a
notion of Finsler connections which named ``$(HM',S,T)$-Cartan
connections''. If, in particular, $HM'$ is the canonical nonlinear
connection $GM'$ of $F^m$ and $S=T=0$, the Finsler connection is
called the Cartan connection and it is denoted by $FC^\ast=(GM',
\nabla^\ast)$ (see \cite{[BeFa]}). They showed that $\nabla^\ast$
is the projection of the Levi-Civita connection of the Sasaki
metric $G$ on the vertical vector bundle. Also they proved that
the associate linear connection $\mathcal{D}^\ast$ to the Cartan
connection $FC^\ast$ is a metric linear connection with respect to
$G$~\cite{[BeFa]}. In this paper we obtain the $(HM',S,T)$-Cartan
connections by using the natural almost complex structure arising
from the nonlinear connection $HM'$, then the natural almost
complex linear connection associated to a $(HM',S,T)$-Cartan
connection is def\/ined. We prove that the natural almost complex
linear connection associated to a $(HM',S,T)$-Cartan connection is
a metric linear connection with respect to the Sasaki metric $G$.
K\"ahler and para-K\"ahler structures associated with Finsler
spaces and their relations with f\/lag curvature were studied by
M.~Crampin and B.Y.~Wu (see \cite{[Crampin], [Wu]}). They have
found some interesting results on this matter. In \cite{[Wu]},
B.Y.~Wu gives some equivalent statements to the K\"ahlerity of
$(M', G, J)$. In the present paper we give other conditions for
the K\"ahlerity of $(M', G, J)$, which extend the previous
results.

\section[The associate almost complex structure to $HM'$]{The associate almost complex structure to $\boldsymbol{HM'}$}

Let $M$ be a real $m$-dimensional smooth manifold and $TM$ be the
tangent bundle of $M$. Let $M'$ be a nonempty open submanifold of
$TM$ such that $\pi(M')=M$ and $\theta(M)\cap M'=\varnothing$,
where $\theta$ is the zero section of $TM$. Suppose that
$F^m=(M,M',F^{\ast})$ is a pseudo-Finsler manifold where
$F^{\ast}: M'\longrightarrow {\Bbb{R}}$ is a smooth function which
in any coordinate system $\{({\mathcal{U}}',\Phi'):x^i,y^i\}$ in
$M'$, the following conditions are fulf\/illed:
\begin{itemize}\itemsep=0pt
    \item $F^\ast$ is positively homogeneous of degree two with
    respect to $(y^1,\dots,y^m)$, i.e., we have
    \[
     F^\ast(x^1,\dots,x^m,ky^1,\dots,ky^m)=k^2F^\ast(x^1,\dots,x^m,y^1,\dots,y^m)
    \]
    for any point $(x,y)\in(\Phi',{\mathcal{U}}')$ and $k>0$.
    \item At any point $(x,y)\in(\Phi',{\mathcal{U}}')$, $g_{ij}$
    are the components of a quadratic form on ${\Bbb{R}}^m$ with
    $q$ negative eigenvalues and $m-q$ positive eigenvalues,
    $0<q<m$ (see \cite{[BeFa]}).
\end{itemize}

Consider the tangent mapping $\pi_\ast: TM'\longrightarrow TM$ of
the submersion $\pi:M'\longrightarrow M$ and def\/ine the vector
bundle $VM'=\ker \pi_\ast$. A complementary distribution $HM'$ to
$VM'$ in $TM'$ is called a nonlinear connection or a horizontal
distribution on $M'$
\[
    TM'=HM'\oplus VM'.
\]
A nonlinear connection $HM'$ enables us to def\/ine an almost
complex structure on $M'$ as follows:
\begin{gather*}
    J:\Gamma(TM')\longrightarrow\Gamma(TM'), \\
   J\left(\frac{\delta}{\delta x^i}\right)=-\frac{\partial}{\partial y^i},\qquad
   J\left(\frac{\partial}{\partial y^i}\right)=\frac{\delta}{\delta x^i},
\end{gather*}
where $\big\{\frac{\delta}{\delta x^i}=\frac{\partial}{\partial
x^i}-N^j_i(x,y)\frac{\partial}{\partial y^j},
\frac{\partial}{\partial y^i}\big\}$ is assumed as a local frame
f\/ield of $TM'$ and $\Gamma(TM')$ is the space of smooth sections
of the vector bundle $TM'$. We call $J$ the associate almost
complex structure to $HM'$. Obviously we have $J^2=-Id_{TM'}$,
also we can assume the conjugate of $J$, $J'=-J$, as an almost
complex structure. Now we give the following proposition which was
proved by B.Y.~Wu~\cite{[Wu]}.
\begin{proposition}\label{proposition1}
Let $F^m=(M, M', F)$ be a Finsler manifold. Then the following
statements are mutually equivalent:
\begin{enumerate}\itemsep=0pt
    \item[{\rm 1)}] $F^m=(M, M', F)$ has zero flag curvature;
    \item[{\rm 2)}] $J$ is integrable;
    \item[{\rm 3)}] $\nabla J=0$, where $\nabla$ is the Levi-Civita connection of the Sasaki
    metric $G$;
    \item[{\rm 4)}] $(M', J, G)$ is K\"ahlerian.
\end{enumerate}
\end{proposition}

\begin{corollary}
Let the associate almost complex structure to $J$ (or $J'$) be a
complex structure; then we have
\begin{gather*}
    \frac{\delta N^k_j}{\delta x^i}=\frac{\delta N^k_i}{\delta x^j},\qquad
    \frac{\partial N^k_j}{\partial y^i}=\frac{\partial N^k_i}{\partial y^j} .
\end{gather*}
So in this case the horizontal distribution is integrable.
\end{corollary}

\section[$(HM',S,T)$-Cartan connection by using the associate almost complex
structure $J$]{$\boldsymbol{(HM',S,T)}$-Cartan connection by using
the associate almost complex structure $\boldsymbol{J}$}

In this section we give another way to def\/ine $(HM',S,T)$-Cartan
connection by using the associate almost complex structure $J$ on
$M'$. Then we study the K\"ahlerity of $(M',J,G)$, where $G$ is
the Sasaki metric and $F^m=(M,M',F)$ is a Finsler manifold.

Let $F^m=(M,M',F^{\ast})$ be a pseudo-Finsler manifold. Then a
Finsler connection on $F^m$ is a~pair $FC=(HM',\nabla)$ where
$HM'$ is a nonlinear connection on $M'$ and $\nabla$ is a linear
connection on the vertical vector bundle $VM'$
(see~\cite{[BeFa]}).
\begin{theorem}\label{theorem1}
Let $\nabla$ be a $FC$ on $M'$. The differential operator
$\mathcal{D}$ defined by
\[
{\mathcal{D}}_XY=\nabla_XvY-J\nabla_XJhY \qquad \forall\;
X,Y\in\Gamma(TM')
\]
is a linear connection on $M'$. Also $J$ is parallel with respect
to ${\mathcal{D}}$, that is
\[
({\mathcal{D}}_XJ)Y=0\qquad \forall \; X,Y\in\Gamma(TM').
\]
We call ${\mathcal{D}}$ the natural almost complex linear
connection associated to $FC$ $\nabla$ on $M'$.
\end{theorem}

\begin{proof}
For any $X,Y,Z\in\Gamma(TM')$ and $f\in {\mathcal{C}}^\infty (M')$
we have
\begin{gather*}
{\mathcal{D}}_{fX+Y}Z=f\nabla_XvZ+\nabla_YvZ-J(f\nabla_XJhZ+\nabla_YJhZ)\nonumber\\
\phantom{{\mathcal{D}}_{fX+Y}Z}{}
=f(\nabla_XvZ-J\nabla_XJhZ)+\nabla_YvZ-J\nabla_YJhZ
=f{\mathcal{D}}_XZ+{\mathcal{D}}_YZ,\\
{\mathcal{D}}_X(fY+Z)=Xf(vY+hY)+f(\nabla_XvY-J\nabla_XJhY)+\nabla_XvZ-J\nabla_XJhZ \\
\phantom{{\mathcal{D}}_X(fY+Z)}{} =
(Xf)Y+f{\mathcal{D}}_XY+{\mathcal{D}}_XZ.
\end{gather*}
Therefore ${\mathcal{D}}$ is a linear connection on $M'$.

Also we have
\begin{gather*}
({\mathcal{D}}_XJ)(Z)={\mathcal{D}}_X(J(Z))-J({\mathcal{D}}_XZ)\\
\phantom{({\mathcal{D}}_XJ)(Z)} {} =\nabla_XvJ(Z)-J\nabla_XJ(h(J(Z)))-J\nabla_XvZ-\nabla_XJhZ\\
\phantom{({\mathcal{D}}_XJ)(Z)} {}
=\nabla_X\left(-Z^i\frac{\partial}{\partial
y^i}\right)-J\nabla_X\left(-\tilde{Z}^i\frac{\partial}{\partial
y^i}\right)-J\nabla_X\left(\tilde{Z}^i\frac{\partial}{\partial
y^i}\right) -\nabla_X\left(-Z^i\frac{\partial}{\partial
y^i}\right) = 0,
\end{gather*}
where in local coordinates $Z=Z^i\frac{\delta}{\delta
x^i}+\tilde{Z}^i\frac{\partial}{\partial y^i}$.
\end{proof}

Note that the torsion of ${\mathcal{D}}$ is given by the following
expression:
\begin{gather}
\label{DTorsion}
T^{\mathcal{D}}(X,Y)=(\nabla_XvY-\nabla_YvX-v[X,Y])
-J(\nabla_XJhY-\nabla_YJhX-Jh[X,Y]).
\end{gather}

\begin{theorem}\label{theorem2}
Let $HM'$ be a nonlinear connection on $M'$ and $S$ and $T$ be any
two skew-symmetric Finsler tensor fields of type $(1,2)$ on $F^m$.
Then there exists a unique linear connection~$\nabla$ on~$VM'$
satisfying the conditions:
\begin{enumerate}\itemsep=0pt
\item[{\rm (i)}] $\nabla$ is a metric connection; \item[{\rm
(ii)}] $T^{\mathcal{D}}$, $S$ and $T$ satisfy
\[
{\rm (a)}\quad T^{\mathcal{D}}(vX,vY)=S(vX,vY), \qquad {\rm
(b)}\quad hT^{\mathcal{D}}(hX,hY)=JT(JhX,JhY)
\]
for any $X,Y\in\Gamma(TM')$, where $J$ is the associate almost
complex structure to $HM'$.
\end{enumerate}
\end{theorem}

\begin{proof}
This proof is similar to \cite{[BeFa]}. We def\/ine a linear
connection $\nabla$ on $VM'$ by using $g$, $h$, $v$, $J$, $S$ and
$T$ in the following way. For any $X, Y, Z \in\Gamma(TM')$ let
\begin{gather}
  2g(\nabla_{vX}vY,vZ)= vX(g(vY,vZ))+vY(g(vZ,vX))-vZ(g(vX,vY))\nonumber\\
  \qquad {} + g(vY,[vZ,vX])+g(vZ,[vX,vY])-g(vX,[vY,vZ])+g(vY,S(vZ,vX))\nonumber\\
  \qquad {} +g(vZ,S(vX,vY))-g(vX,S(vY,vZ))\label{eq1}
\end{gather}
and
\begin{gather}
   2g(\nabla_{hX}JhY, JhZ)= hX(g(JhY,JhZ))+hY(g(JhZ,JhX))\nonumber\\
   \qquad {}  - hZ(g(JhX,JhY))+g(JhY,Jh[hZ,hX])+ g(JhZ,Jh[hX,hY]) \nonumber\\
   \qquad {}   -g(JhX,Jh[hY,hZ])+g(JhY,T(JhZ,JhX))\nonumber  \\
   \qquad {}  +g(JhZ,T(JhX,JhY))-g(JhX,T(JhY,JhZ)).\label{eq2}
\end{gather}
Then for any $X, Y, Z \in\Gamma(TM')$ we have
\begin{gather}
  (\nabla_Xg)(vY,vZ) = (\nabla_{vX+hX}g)(vY,vZ)\nonumber\\
                     \qquad {}  = vX(g(vY,vZ))-g(\nabla_{vX}vY,vZ)-g(vY,\nabla_{vX}vZ)+hX(g(vY,vZ))\nonumber\\
                     \qquad {}   -g(\nabla_{hX}vY,vZ)-g(vY,\nabla_{hX}vZ)=0.\nonumber
\end{gather}
The above computation shows that the connection $\nabla$ def\/ined
by~\eqref{eq1} and~\eqref{eq2} is a metric connection.

Locally we set $\nabla_{\frac{\delta}{\delta
x^j}}\frac{\partial}{\partial
y^i}=F^k_{ij}(x,y)\frac{\partial}{\partial y^k}$,
$\nabla_{\frac{\partial}{\partial y^j}}\frac{\partial}{\partial
y^i}=C^k_{ij}(x,y)\frac{\partial}{\partial y^k}$,
$S(\frac{\partial}{\partial y^j}, \frac{\partial}{\partial
y^i})=S^k_{ij}\frac{\partial}{\partial y^k}$ and
$T(\frac{\partial}{\partial y^j}, \frac{\partial}{\partial
y^i})=T^k_{ij}\frac{\partial}{\partial y^k}$.

Now in~\eqref{eq1} let $X=\frac{\partial}{\partial y^j}$,
$Y=\frac{\partial}{\partial y^i}$ and $Z=\frac{\partial}{\partial
y^l}$. After performing some computations we obtain the following
expression for the coef\/f\/icients $C^m_{ij}$:
\begin{gather*}
  C^m_{ij}=\frac{1}{2}\left\{\frac{\partial g_{il}}{\partial y^j}
  +\frac{\partial g_{lj}}{\partial y^i}-\frac{\partial g_{ji}}{\partial y^l}
  +S^h_{jl}g_{ih}+S^h_{ij}g_{lh}-S^h_{li}g_{jh}\right\}g^{lm}.
\end{gather*}
Also in~\eqref{eq2} let $X=\frac{\delta}{\delta x^j}$,
$Y=\frac{\delta}{\delta x^i}$ and $Z=\frac{\delta}{\delta x^l}$.
Then we can obtain the following expression for the coef\/f\/icients
$F^m_{ij}$:
\begin{gather*}
  F^m_{ij}=\frac{1}{2}\left\{\frac{\delta g_{il}}{\delta x^j}
  +\frac{\delta g_{lj}}{\delta x^i}-\frac{\delta g_{ji}}{\delta x^l}
  -T^h_{jl}g_{ih}-T^h_{ij}g_{lh}+T^h_{li}g_{jh}\right\}g^{lm}.
\end{gather*}
By using the relations $J\circ v=h\circ J$, $v\circ J=J\circ h$
and~\eqref{DTorsion} we have
\begin{gather}
T^{\mathcal{D}}(vX,vY)=\nabla_{vX}vY-\nabla_{vY}vX-[vX,vY],\label{vDTorsion}\\
hT^{\mathcal{D}}(hX,hY)=J(\nabla_{hY}JhX-\nabla_{hX}JhY+Jh[hX,hY])\label{hDTorsion}.
\end{gather}
Suppose that $X,Y\in\Gamma(TM')$ are two arbitrary vector f\/ields
on $M'$. In local coordinates, let $X=X^i\frac{\delta}{\delta
x^i}+\tilde{X}^i\frac{\partial}{\partial y^i}$ and
$Y=Y^i\frac{\delta}{\delta
x^i}+\tilde{Y}^i\frac{\partial}{\partial y^i}$, after performing
some computations we have:
\begin{gather}
  T^{\mathcal{D}}\left(\tilde{X}^i\frac{\partial}{\partial y^i}
  ,\tilde{Y}^i\frac{\partial}{\partial y^i}\right)=S\left(\tilde{X}^i\frac{\partial}{\partial y^i}
  ,\tilde{Y}^i\frac{\partial}{\partial y^i}\right),\label{eq5}\\
  hT^{\mathcal{D}}\left(X^i\frac{\delta}{\delta
  x^i},Y^i\frac{\delta}{\delta x^i}\right)=JT\left(J\left(X^i\frac{\delta}{\delta
  x^i}\right),J\left(Y^i\frac{\delta}{\delta x^i}\right)\right)\label{eq6}.
\end{gather}
The relations~\eqref{eq5} and~\eqref{eq6} show that $\nabla$
satisf\/ies (ii) of Theorem~\ref{theorem2}.

Now let $\tilde{\nabla}$ be another linear connection on $VM'$
which satisf\/ies (i) and (ii). By using the relations (i), (ii),
\eqref{vDTorsion} and~\eqref{hDTorsion} for $\tilde{\nabla}$ we
have the following expressions:
\begin{gather}
  vX(g(vY,vZ))+vY(g(vZ,vX))-vZ(g(vX,vY))\nonumber\\
  \qquad {} =g(\tilde{\nabla}_{vX}vY+\tilde{\nabla}_{vX}vY-T^{\mathcal{D}}(vX,vY)-[vX,vY],vZ)\nonumber\\
  \qquad {} +g(T^{\mathcal{D}}(vX,vZ)+[vX,vZ],vY) +g(T^{\mathcal{D}}(vY,vZ)+[vY,vZ],vX),\label{eq7}\\
  hX(g(vJY,vJZ))+hY(g(vJZ,vJX))-hZ(g(vJX,vJY))\nonumber\\
  \qquad {} =g(\tilde{\nabla}_{hX}JhY+\tilde{\nabla}_{hX}JhY-JT(JhX,JhY)-Jh[hX,hY],JhZ)\label{eq8}\\
  \qquad {} +g(JT(JhX,JhZ)+Jh[hX,hZ],JhY) +g(JT(JhY,JhZ)+Jh[hY,hZ],JhX).\nonumber
\end{gather}
The relations~\eqref{eq7} and~\eqref{eq8} show that
$\tilde{\nabla}$ satisf\/ies~\eqref{eq1} and~\eqref{eq2},
respectively. Therefore $\nabla=\tilde{\nabla}$.
\end{proof}

The Finsler connection $FC=(HM',\nabla)$ where $\nabla$ is given
by Theorem~\ref{theorem2} is called the $(HM',S,T)$-Cartan
connection (see \cite{[BeFa], [BeFa1]}) which in this case is
obtained by the associate almost complex structure to $HM'$. If,
in particular, $HM'$ is just the canonical nonlinear
connection~$GM'$ of ${\Bbb{F}}^m$ (for more details about $GM'$
see~\cite{[BeFa]}) and $S=T=0$, the $FC$ is called the Cartan
connection and it is denoted by $FC^\ast=(GM',\nabla^\ast)$.

By means of the pseudo-Riemannian metric $g$ on $VM'$  we consider
a pseudo-Riemannian metric on the vector bundle $TM'$ similar to
the Sasaki one and denote it by $G$, that is
\[
  G=g_{ij}(x,y)dx^idx^j+g_{ij}(x,y)\delta y^i\delta y^j,
\]
where $\delta y^i=dy^i+N^i_j(x,y)dx^j$. Denote by $\nabla'$ the
Levi-Civita connection on $M'$ with respect to~$G$. A.~Bejancu and
H.R.~Farran showed $\nabla^\ast$ is the projection of the
Levi-Civita connection~$\nabla'$ on the vertical vector bundle
also they proved the following theorem (see~\cite{[BeFa]}).

\begin{theorem}\label{theorem3}
The associate linear connection ${\mathcal{D}}^\ast$ to the Cartan
connection $FC^\ast=(GM',\nabla^\ast)$ is a metric linear
connection with respect to $G$.
\end{theorem}

Now we give the following theorem which shows the natural almost
complex linear connections associated to $(HM',S,T)$-Cartan
connections are metric linear connections with respect to $G$.

\begin{theorem}\label{theorem4}
The natural almost complex linear connection ${\mathcal{D}}$
associated to a $(HM',S,T)$-Cartan connection $FC=(HM',\nabla)$ is
a metric linear connection with respect to $G$.
\end{theorem}

\begin{proof}
For any $X,X_1,X_2\in\Gamma(TM')$ we have
\begin{gather}
  {\mathcal{D}}_XG(X_1,X_2)= XG(X_1,X_2)-G({\mathcal{D}}_XX_1,X_2)-G(X_1,{\mathcal{D}}_XX_2)\nonumber\\
 \phantom{{\mathcal{D}}_XG(X_1,X_2)}{}  =X(G(X_1,X_2))-G(\nabla_XvX_1,X_2)+G(J\nabla_XJhX_1,X_2) \nonumber\\
\phantom{{\mathcal{D}}_XG(X_1,X_2)=}{}
-G(X_1,\nabla_XvX_2)+G(X_1,J\nabla_XJhX_2).\label{eq9}
\end{gather}
By~\eqref{eq9} and this fact that $S$ and $T$ are skew-symmetric
we have:
\begin{gather*}
    {\mathcal{D}}_{\frac{\partial}{\partial y^i}}G\left(\frac{\partial}{\partial y^j},\frac{\delta}{\delta x^k}\right)=
    {\mathcal{D}}_{\frac{\delta}{\delta x^i}}G\left(\frac{\partial}{\partial y^j},\frac{\delta}{\delta x^k}\right)=0,\\
    {\mathcal{D}}_{\frac{\partial}{\partial y^i}}G\left(\frac{\partial}{\partial y^j},\frac{\partial}{\partial y^k}\right)=
    {\mathcal{D}}_{\frac{\partial}{\partial y^i}}G\left(\frac{\delta}{\delta x^j},\frac{\delta}{\delta x^k}\right)
    =\frac{\partial g_{jk}}{\partial y^i}-C^h_{ji}g_{hk}-C^h_{ki}g_{jh}=0 ,\\
    {\mathcal{D}}_{\frac{\delta}{\delta x^i}}G\left(\frac{\partial}{\partial y^j},\frac{\partial}{\partial y^k}\right)
    = {\mathcal{D}}_{\frac{\delta}{\delta x^i}}G\left(\frac{\delta}{\delta x^j},\frac{\delta}{\delta x^k}\right)
    =\frac{\delta g_{jk}}{\delta x^i}-F^h_{ji}g_{hk}-F^h_{ki}g_{jh}
    =0.
\end{gather*}
Therefore ${\mathcal{D}}_X G=0$ for any $X\in \Gamma(TM')$.
\end{proof}

Let $F^m=(M,M',F)$ be a Finsler manifold. We can easily check that
the pair $(J,G)$ def\/ines an almost Hermitian metric on $M'$. In
the following theorem we give a suf\/f\/icient condition for Finsler
tensor f\/ields $S$ and $T$ such that ${\mathcal{D}}$ be the
Levi-Civita connection arising from $G$.

\begin{theorem}\label{theorem5}
The natural almost complex linear connection ${\mathcal{D}}$
associated to a $(HM',S,T)$-Cartan connection $FC=(HM',\nabla)$ is
the Levi-Civita connection arising from $G$ if
$T^{\mathcal{D}}(X,Y)=0$ for any $X,Y\in\Gamma(TM')$ or
equivalently if
\begin{gather*}
  S=T=0,\qquad   C^k_{ij}=R^k_{ij}=0,\qquad F^k_{ij}=\frac{\partial N^k_j}{\partial
  y^i},
\end{gather*}
where $R^k_{ij}=\frac{\delta N^k_i}{\delta x^j}-\frac{\delta
N^k_j}{\delta x^i}$.
\end{theorem}

\begin{proof}
By Theorem~\ref{theorem4}, ${\mathcal{D}}$ is a metric linear
connection with respect to $G$. Therefore if $T^{\mathcal{D}}=0$
then ${\mathcal{D}}$ is the Levi-Civita connection. In local
coordinates we have
\begin{gather*}
  T^{\mathcal{D}}\left(\frac{\partial}{\partial y^j},\frac{\partial}{\partial y^i}\right)
  = S^k_{ij}\frac{\partial}{\partial y^k},\\
  T^{\mathcal{D}}\left(\frac{\partial}{\partial y^i},\frac{\delta}{\delta x^j}\right)= C^k_{ji}\frac{\delta}{\delta
  x^k}+\left(\frac{\partial N^k_j}{\partial y^i}-F^k_{ij}\right)\frac{\partial}{\partial
  y^k},\\
  T^{\mathcal{D}}\left(\frac{\delta}{\delta x^i},\frac{\delta}{\delta x^j}\right)= T^k_{ij}\frac{\delta}{\delta
  x^k}+\left(\frac{\delta N^k_j}{\delta x^i}-\frac{\delta N^k_i}{\delta x^j}\right)\frac{\partial}{\partial
  y^k}.
\end{gather*}
Therefore the proof is completed.
\end{proof}

\begin{corollary}\label{corollary2}
If $T^{\mathcal{D}}=0$ then $(M',J,G)$ is a K\"ahler manifold.
\end{corollary}

\begin{proof}
If $T^{\mathcal{D}}=0$ then ${\mathcal{D}}$ is the Levi-Civita
connection of $G$. Also $J$ is parallel with respect to
${\mathcal{D}}$. Therefore ${\mathcal{D}}$ (the Levi-Civita
connection of $G$) is almost complex. Consequently by using
Theorem 4.3 of \cite{[KoNo]}, $(M',J,G)$ is a K\"ahler manifold.
\end{proof}

We know that the almost Hermitian manifold $(M',J,G)$ is an almost
K\"ahler manifold if and only if the fundamental $2$-form $\Phi$
is closed ($\Phi$ is def\/ined by $\Phi(X,Y)=G(X,JY)$ for all
$X,Y\in \Gamma(TM')$). Therefore we can give the following
theorem.

\begin{theorem}\label{theorem6}
The almost Hermitian manifold $(M',J,G)$ is an almost K\"ahler
manifold if and only if
\begin{gather}\label{eq10}
  \frac{\delta g_{ik}}{\delta x^j}+\frac{\partial N^h_k}{\partial y^i}g_{hj}-
  \left(\frac{\delta g_{ij}}{\delta x^k}+\frac{\partial N^h_j}{\partial y^i}g_{hk} \right)=0
\end{gather}
and
\begin{gather}\label{eq11}
  R^h_{ij}g_{hk}-R^h_{ik}g_{hj}+R^h_{jk}g_{hi}=0.
\end{gather}
\end{theorem}

\begin{proof}
Let $X_0, X_1, X_2\in\Gamma(TM')$. Then we have
\begin{gather*}
  d\Phi(X_0,X_1,X_2) = X_0G(X_1,JX_2)-X_1G(X_0,JX_2)+X_2G(X_0,JX_1) \\
  \phantom{d\Phi(X_0,X_1,X_2) =} {} -G([X_0,X_1],JX_2)+G([X_0,X_2],JX_1)-G([X_1,X_2],JX_0).
\end{gather*}
By using the above relation in local coordinates we have:
\begin{gather*}
  d\Phi\left(\frac{\partial}{\partial y^i},\frac{\partial}{\partial y^j},\frac{\partial}{\partial y^k}\right)
  = d\Phi\left(\frac{\partial}{\partial y^i},\frac{\partial}{\partial y^j},\frac{\delta}{\delta x^k}\right)=0,\\
  d\Phi\left(\frac{\partial}{\partial y^i},\frac{\delta}{\delta x^j},\frac{\delta}{\delta
  x^k}\right)=\frac{\delta g_{ik}}{\delta x^j}+\frac{\partial N^h_k}{\partial y^i}g_{hj}-
  \left(\frac{\delta g_{ij}}{\delta x^k}+\frac{\partial N^h_j}{\partial y^i}g_{hk}\right),\\
  d\Phi\left(\frac{\delta}{\delta x^i},\frac{\delta}{\delta x^j},\frac{\delta}{\delta
  x^k}\right)=\left(\frac{\delta N^h_i}{\delta x^j}-\frac{\delta N^h_j}{\delta
  x^i}\right)g_{hk}-\left(\frac{\delta N^h_i}{\delta x^k}-\frac{\delta N^h_k}{\delta
  x^i}\right)g_{hj}+\left(\frac{\delta N^h_j}{\delta x^k}-\frac{\delta N^h_k}{\delta
  x^j}\right)g_{hi}.
\end{gather*}
Therefore the fundamental $2$-form $\Phi$ is closed if and only if
the equations~\eqref{eq10} and~\eqref{eq11} are conf\/irmed.
\end{proof}

Now, by using Proposition~\ref{proposition1} and
Corollary~\ref{corollary2}, we have the following corollary.
\begin{corollary}
Let $F^m=(M, M', F)$ be a Finsler manifold. If $T^{\mathcal{D}}=0$
then,
\begin{enumerate}\itemsep=0pt
    \item[{\rm 1)}] $F^m=(M, M', F)$ has zero flag curvature;
    \item[{\rm 2)}] $J$ is integrable;
    \item[{\rm 3)}] $\nabla J=0$, where $\nabla$ is the Levi-Civita connection of the Sasaki
    metric $G$;
    \item[{\rm 4)}] $(M', J, G)$ is K\"ahlerian.
\end{enumerate}

\end{corollary}

\subsection*{Acknowledgements}

The authors are grateful to the referees for their valuable
suggestions on this paper.

\LastPageEnding

\end{document}